\documentclass[a4paper,12pt]{article}

\usepackage{amsmath, amsthm, amssymb}
\usepackage{enumerate}
\usepackage{graphicx}
\usepackage{natbib}
\usepackage[top=1in, bottom=1in, left=1in, right=1in]{geometry}

\title{Law of large numbers unifying Maxwell-Boltzmann, Bose-Einstein and Zipf-Mandelbort distributions, and related fluctuations.}
\author{Tomasz M. \L api\'nski \thanks{Faculty of Applied Physics and Mathematics, Gda\'nsk University of Technology, ul. Gabriela Narutowicza 11/12, 80-233 Gda\'nsk, Poland, email: 84tomek@gmail.com}}

\begin{document}

\newtheorem{remark}{Remark}
\newtheorem{preposition}{Preposition}
\newtheorem{definition}{Definition}
\newtheorem{proposition}{Proposition}
\newtheorem{theorem}{Theorem}
\newtheorem{lemma}{Lemma}
\newtheorem{proofoflemma}{Proof of Lemma}
\newtheorem{proofofproposition}{Proof of Proposition}
\newtheorem{proofoflpreposition}{Proof of Preposition}

\maketitle

\begin{abstract}
	We consider a system composed of a fixed number of particles with total energy smaller or equal to some prescribed value. The particles are non-interacting, indistinguishable and distributed over fixed number of energy levels. The energy levels are degenerate and degeneracy is a function of the number of particles. Three cases of the degeneracy function is considered. It can increase with either the same rate as the number of particles or slower, or faster. We find useful properties of the entropy of the system and solve related entropy optimization problem. It turned out, there are several solutions. Depending on the magnitude of total energy, the maximum of the entropy can be in the interior of system's state space or on the boundary. On the boundary it can have further three cases depending on the degeneracy function. The main result, Law of Large Numbers yields the most probable state of the system, which equals to the point of maximum of the entropy. This point can be either Maxwell-Boltzmann statistics or Bose-Einstein statistics, or Zipf-Mandelbort law. We also find the limiting laws for the fluctuations. These laws are different for different cases of the entropy's maximum. They can be mixture of Normal, Exponential and Discrete distributions. Explicit rates of convergence of moment generating functions are provided for all the theorems. The overview of possible applications is provided in the last section.
\end{abstract}


\begin{msc}
	60F05, 82B03
\end{msc}
\\
\\
\copyright 2021. This manuscript version is made available under the CC-BY-NC-ND 4.0 license http://creativecommons.org/licenses/by-nc-nd/4.0/.\\
\\
https://doi.org/10.1016/j.physa.2021.125909
\\


\maketitle


\section{Introduction}
	\qquad The system under consideration is composed of $N$ non-interacting, indistinguishable particles with total energy smaller or equal to EN. The constant $E$ is a maximal average energy per particle. The particles are distributed over $m$ energy levels and those levels have energies $\varepsilon=(\varepsilon_{1},\varepsilon_{2},\ldots,\varepsilon_{m})$ where $\varepsilon_{1}<\varepsilon_{2}<\ldots<\varepsilon_{m}$. A single system's state is represented by the vector of level's occupation numbers, that is, $(N_{1},N_{2},\ldots,N_{m})$. Since the number of particles is fixed to $N$ and energy is smaller or equal to $EN$, the occupation numbers satisfy the following constrains
	\begin{align}
		\label{particle_system_ensamble_constraint_1}
		&N=N_{1}+N_{2}+\ldots+N_{m},\\
		\label{particle_system_ensamble_constraint_2}
		&EN \geq \varepsilon_{1}N_{1}+\varepsilon_{2}N_{2}+\ldots+\varepsilon_{m}N_{m}.
	\end{align} 
	Each energy level has $G_{i}$ degeneracy, where $i=1,\ldots,m$. Hence, the entropy of the particular system's state is the following
	\begin{equation}
		\label{introduction_entropy_standard}
		S(N_{1},\ldots,N_{m})=\sum_{i=1}^{m}\ln\frac{(N_{i}+G_{i}-1)!}{N_{i}!(G_{i}-1)!}.
	\end{equation}
	The total degeneracy of the levels is equal to $G=\sum_{i=1}^{m}G_{i}$. We assume total degeneracy is an increasing function of the number of particles. We study three regimes of the behavior of $G=G(N)$ 
	\begin{align}
		\label{introduction_G_cases}
		1)\ &\lim_{N\to\infty}\frac{G(N)}{N}=\infty,\notag\\
		2)\ &\frac{G(N)}{N}=\alpha+\beta(N),\ \text{where}\ \alpha>0, \lim_{N\to\infty}\beta(N)=0,\\
		3)\ &\lim_{N\to\infty}\frac{G(N)}{N}=0.\notag
	\end{align}
	Moreover, for each $N$, the components $G_{i}$ are equally weighted. That is, for all $N$, $G_{i}=g_{i}G(N)$, for $i=1,\ldots,m$ and some constants $g_{i}$ such that $\sum_{i=1}^{m}g_{i}=1$.\\
	\indent Next, for each system's state let us introduce a vector of weights of the occupation numbers, that is
	\begin{equation}
		\label{introduction_vector_of_weights}
		x=(x_{1},x_{2},\ldots,x_{m})=\bigg(\frac{N_{1}}{N},\frac{N_{2}}{N},\ldots,\frac{N_{2}}{N}\bigg).
	\end{equation}
	In this context, the entropy (\ref{introduction_entropy_standard}) is equal
	\begin{equation}
		\label{introduction_entropy}
		S(x,N)=\sum_{i=1}^{m}\ln\frac{(x_{i}N+g_{i}G(N)-1)!}{(x_{i}N)!(g_{i}G(N)-1)!}.
	\end{equation}
	Now, let us define a set $\mathcal{A}_{E}$ composed of the vectors $x\in\mathbb{R}^{m}$, such that the following constraints are valid
	\begin{align}
		\label{particle_system_domain_constraint_1}
		&1=x_{1}+x_{2}+\ldots+x_{m},\\
		\label{particle_system_domain_constraint_2}
		&E \geq \varepsilon_{1}x_{1}+\varepsilon_{2}x_{2}+\ldots+\varepsilon_{m}x_{m},\\
		&x_{i}\geq 0, \quad i=1,\ldots,m, \notag
	\end{align}
	and also define a set $L_{N}:=\{\frac{\mathcal{N}}{N},\mathcal{N}\in\mathbb{N}^{m}\}$ for any $N\geq N_{0}$, $N_{0}\in\mathbb{Z}_{+}$. We assume that $0\in\mathbb{N}^{m}$. Taking into account constraints (\ref{particle_system_ensamble_constraint_1}), (\ref{particle_system_ensamble_constraint_2}) and representation of single state (\ref{introduction_vector_of_weights}), it is clear that the state space of the system is equal to the set $\mathcal{A}_{E}\cap L_{N}$. Moreover, it can be easily verified that  when $E=\varepsilon_{1}$, the set $\mathcal{A}_{E}$ is a single point. When $E<\varepsilon_{1}$, $\mathcal{A}_{E}$ is an empty set. So, let us assume the parameter $E$ is larger than $\varepsilon_{1}$.\\
	\indent Our first result is an approximation of the entropy (\ref{introduction_entropy}) by the expression suitable for further development. We also identify several useful properties of the resulting approximation and the entropy itself.\\
	\indent Then, the next result is about entropy's point of maximum for all possible regimes of $E$ and all the cases of the function $G(N)$.\\
	\indent Then we consider a discrete random vector $X(N):=\big(X_{1}(N),X_{2}(N),\ldots,X_{m}(N)\big)$ where $X(N)\in\mathcal{A}_{E}\cap L_{N}$. We use fundamental postulate of statistical mechanics that system's microstates are equally probable, see e.g. \cite{statistical_physics_pathria_book} and \cite{statistical_physics_reif_book}, in order to define the pmf of $X(N)$ as
	\begin{equation}
		\label{particle_system_pmf}
		P(X(N)=x):=\frac{e^{S(x,N)}}{\sum_{\mathcal{A}_{E}\cap L_{N}}e^{S(y,N)}},
	\end{equation}
	where $x$ is a single system's state and the entropy $S(x,N)$ is given by (\ref{introduction_entropy}).\\
	\indent The main result of this work is the law of large numbers that yields the most probable state of the random vector $X(N)$ as the number of particles is large enough. This state is the point of maximum of the entropy (\ref{introduction_entropy}). When this point of maximum is in the interior of the domain, it is equal $x^{*}=g$, where $g=(g_{1},g_{2},\ldots,g_{m})$. When the maximum is on the boundary, we obtain different points of maximum for each case of $G(N)$. These points are distributions: Maxwell-Boltzmann, Bose-Einstein and Zipf-Mandelbort, respectively for three cases in (\ref{introduction_G_cases}). Our second result yields the distributions of fluctuations from the most probable system's state. They are also different for two cases of entropy's maximum. When the maximum is in the interior of the domain, fluctuations have Normal distribution. When the maximum is on the boundary, there can be further two cases depending on the degeneracy function. For the first and the second case in (\ref{introduction_G_cases}), the fluctuations distribution is Exponential in the direction orthogonal to the boundary on state space and Normal in other directions. For the third case fluctuations distribution is Discrete in direction orthogonal to the boundary and Normal in other directions. Explicit rates of convergence of moment generating functions are provided for all the limit theorems. \\
	\indent The results presented in this paper are continuation of the work done as a part of the author's PhD thesis \cite{phd}. There the considered system is introduced and the entropy properties are partially developed. Most notably, the thesis lacks a rigorous proof of the limit theorems.\\
	\indent The same system, also with degeneracy depending on the number of particles, was introduced in \cite{maslov_general_theorem}. There the author stated existence of the two cases of maximum of the entropy depending on the total energy. The proof was not provided. For the case when the entropy maximum is on the boundary and for the single case of $G(N)$, the optimization problem similar to ours was solved in \cite{maslov_nonlinear_axioms}. In \cite{maslov_general_theorem} and \cite{maslov_nonlinear_averages} the authors also proved convergence to the three distributions: Maxwell-Boltzmann, Bose-Einstein and Zipf-Mandelbort. However, only for the situation when maximum is on the boundary. Our results on the fluctuations are completely new. The proofs there are lacking rigor and are completely different from ours. For other results on the fluctuations of a related Bose-Einstein condensate see \cite{bose-einstein_condensate}.\\
	\indent Proofs of our limit theorems are based on the results from \cite{lapinski}, which are summarized in the Appendix. It should be mentioned that the results in \cite{Laplace_method_approach_kolokoltsov} and \cite{vassili_notes} significantly influenced the theorems recalled in the Appendix.\\


\section{Entropy properties}
	For some $a>0$, let us define a set $\mathcal{A}_{E}^{(a)}:=\{x:x\in\mathcal{A}_{E}, x_{i}>a; i=1,\ldots,m\}$. 
	\begin{proposition}
		\label{entropy_properties_preposition_approximation}
		The entropy $S(x,N)$ given by (\ref{introduction_entropy}) is three times differentiable and strictly concave function on $\mathcal{A}_{E}$. Moreover, for large enough $N$, and any $x\in\mathcal{A}^{(a)}_{E}$ the following approximation holds
		\begin{equation}
			\label{entropy_properties_preposition_approximation_entropy}
			S(x,N)=h(N)f(x,N)+C(N),\ \text{with}\ f(x,N)=f(x)+\sigma(x,N)\epsilon(N),
		\end{equation}
		where
		\begin{enumerate}[i)]
			\item $f(x,N)$ and $f(x)$ are strictly concave and three times differentiable functions,
			\item $h(N)$ and $C(N)$ are some functions of $N$,
			\item $\sigma(x,N)$ and its derivatives up to third order are uniformly bounded,
		\end{enumerate}
		and for each case of $G(N)$ given by (\ref{introduction_G_cases}) we have the following explicit form of the above functions
		\begin{enumerate}[1)]
			\item 
				\begin{align*}
					&h(N)=N,\ C(N)=N\bigg(\ln\frac{G(N)}{N}+1\bigg)-\frac{1}{2}\ln N,&\\
					&f(x)=\sum_{i=1}^{m}x_{i}\ln\frac{g_{i}}{x_{i}},& \\
					&\epsilon(N)=
					\begin{cases}
						\frac{1}{N}& \qquad \text{if}\ \lim_{N\to\infty}\frac{N^{2}}{G(N)}=0,\\
					 	\frac{N}{G(N)}& \qquad \text{otherwise}.
					\end{cases}
				\end{align*}
			\item 
				\begin{align*}
					&h(N)=N,\ C(N)=-N\sum_{i=1}^{m}g_{i}\alpha\ln g_{i}\alpha-\frac{1}{2}\ln N,&\\
					&f(x)=\sum_{i=1}^{m}\big((x_{i}+g_{i}\alpha)\ln(x_{i}+g_{i}\alpha)-x_{i}\ln x_{i}\big),&\\
					&\epsilon(N)=
					\begin{cases}
						\frac{1}{N}& \qquad \text{if}\ \lim_{N\to\infty}\beta(N)N=0,\\
					 	|\beta(N)|& \qquad \text{otherwise}.
					\end{cases}
				\end{align*}
			\item
				\begin{align*}
					&h(N)=G(N),\ C(N)=G(N)(\ln N+1)-\sum_{i=1}^{m}\bigg(g_{i}G(N)+\frac{1}{2}\bigg)\ln g_{i}G(N)-\ln N,&\\
					&f(x)=\sum_{i=1}^{m}g_{i}\ln x_{i},&\\
					&\epsilon(N)=
					\begin{cases}
						\frac{1}{G(N)}& \qquad \text{if}\ \lim_{N\to\infty}\frac{G(N)^{2}}{N}=0,\\
					 	\frac{G(N)}{N}& \qquad \text{otherwise}.
					\end{cases}
				\end{align*}
		\end{enumerate}
		
	\end{proposition}

	\begin{proof}
		Let us recall Stirling's approximation
		\begin{equation}
			\label{entropy_properties_preposition_approximation_proof_stirling}
			\ln(\lambda !)=\bigg(\lambda+\frac{1}{2}\bigg)\ln\lambda-\lambda+\ln\sqrt{2\pi}+\frac{\omega(\lambda)}{\lambda}.
		\end{equation}
		The remainder $\omega(\lambda)$ is an analytic function and $\omega(\lambda)=O(1)$ as $\lambda\to\infty$.\\
		We use approximation (\ref{entropy_properties_preposition_approximation_proof_stirling}) to obtain
		\begin{align}
			\label{entropy_properties_preposition_approximation_proof_approximation}
			S(x,N)=&\sum_{i=1}^{m}\Bigg[\bigg(x_{i}N+g_{i}G(N)-\frac{1}{2}\bigg)\ln(x_{i}N+g_{i}G(N))-\bigg(x_{i}N+\frac{1}{2}\bigg)\ln x_{i}N\\
			&-\bigg(g_{i}G(N)-\frac{1}{2}\bigg)\ln (g_{i}G(N))+\omega(x_{i},N)\Bigg]\notag,
		\end{align}
		where the remainder term is equal
		\begin{align}
			\label{entropy_properties_preposition_approximation_proof_omega_i}
			\omega(x_{i},&N)=\bigg(x_{i}N+g_{i}G(N)-\frac{1}{2}\bigg)\ln\bigg(1-\frac{1}{x_{i}N+g_{i}G(N)}\bigg)-\ln\sqrt{2\pi}\\
			&-\bigg(g_{i}G(N)-\frac{1}{2}\bigg)\ln\bigg(1-\frac{1}{g_{i}G(N)}\bigg)+\frac{\omega_{1}(x_{1},N)}{x_{i}N+g_{i}G(N)}-\frac{\omega_{2}(x_{i},N)}{x_{i}N}-\frac{\omega_{3}(N)}{g_{i}G(N)-1}.\notag
		\end{align}
		where $\omega_{1}$, $\omega_{2}$, $\omega_{3}$ are the remainders of applied formula (\ref{entropy_properties_preposition_approximation_proof_stirling}). Moreover, for any $a\in(0,1)$, there exists $N_{0}\in\mathbb{Z}_{+}$ such that approximation (\ref{entropy_properties_preposition_approximation_proof_omega_i}) holds for all $x_{i}\in(a,1]$ and all $N\geq N_{0}$. Therefore, for large enough $N$, approximation (\ref{entropy_properties_preposition_approximation_proof_approximation}) holds for any $x\in\mathcal{A}_{E}^{(a)}$.\\
		\indent Let us perform further calculations separately for each case of $G(N)$ given by (\ref{introduction_G_cases})
		\begin{enumerate}[1)]
		\item Using (\ref{entropy_properties_preposition_approximation_proof_approximation}) the function $S(x,N)$ can be presented as
			\begin{equation*}
				S(x,N)=Nf(x,N)+N\bigg(\ln\frac{G(N)}{N}+1\bigg)-\frac{1}{2}\ln N,
			\end{equation*}
				where
			\begin{align*}
				f(x,N)=&\sum_{i=1}^{m}\bigg[x_{i}\ln\frac{g_{i}}{x_{i}}+\bigg(x_{i}+\frac{g_{i}G(N)}{N}-\frac{1}{2N}\bigg)\ln\bigg(1+\frac{x_{i}N}{g_{i}G(N)}\bigg)-x_{i}-\frac{1}{2N}\ln x_{i}\\
				&+\frac{\omega(x_{i},N)}{N
}\bigg].			
			\end{align*}
			Then the function $f(x,N)$ can be presented as
			\begin{equation*}
				f(x,N)=f(x)+\sigma(x,N)\epsilon(N),
			\end{equation*}
			where
			\begin{align*}
				f(x)=&\sum_{i=1}^{m}x_{i}\ln\frac{g_{i}}{x_{i}},\\
				\sigma(x,N)=&\frac{1}{\epsilon(N)}\sum_{i=1}^{m}\Bigg[\bigg(x_{i}+\frac{g_{i}G(N)}{N}-\frac{1}{2N}\bigg)\ln\bigg(1+\frac{x_{i}N}{g_{i}G(N)}\bigg)-x_{i}-\frac{1}{2N}\ln x_{i}+\frac{\omega(x_{i},N)}{N}\Bigg],\\
				\epsilon(N)=&
				\begin{cases}
					\frac{1}{N} \qquad &\text{if}\ \lim_{N\to\infty}\frac{N^{2}}{G(N)}=0,\\
					\frac{N}{G(N)} \qquad &\text{otherwise}.
				\end{cases}
			\end{align*}
		\item Analogously to 1), by (\ref{entropy_properties_preposition_approximation_proof_approximation}),
			\begin{equation*}
				S(x,N)=Nf(x,N)-N\sum_{i=1}^{m}g_{i}\alpha\ln g_{i}\alpha-\frac{1}{2}\ln N,
			\end{equation*}
				where
			\begin{align*}
				f(x,N)=&\sum_{i=1}^{m}\Bigg[\bigg(x_{i}+\frac{g_{i}G(N)}{N}+\frac{1}{2N}\bigg)\ln\bigg(x_{i}+\frac{g_{i}G(N)}{N}\bigg)-x_{i}\ln x_{i}-\frac{1}{2N}\ln x_{i}-\\
				&-\bigg(\frac{g_{i}G(N)}{N}-\frac{1}{2N}\bigg)\ln\frac{g_{i}G(N)}{N}+g_{i}\alpha\ln g_{i}\alpha+\frac{\omega(x_{i},N)}{N}\Bigg],
			\end{align*}
			and also
			\begin{equation*}
				f(x,N)=f(x)+\sigma(x,N)\epsilon(N),
			\end{equation*}
			where
			\begin{align*}
				f(x)=&\sum_{i=1}^{m}(x_{i}+g_{i}\alpha)\ln(x_{i}+g_{i}\alpha)-x_{i}\ln x_{i},\\
				\sigma(x,N)=&\frac{1}{\epsilon(N)}\sum_{i=1}^{m}\Bigg[(x_{i}+g_{i}\alpha)\ln\bigg(1+\frac{g_{i}\beta(N)}{x_{i}+g_{i}\alpha}\bigg)+\bigg(g_{i}\beta(N)+\frac{1}{2N}\bigg)\ln\bigg(x_{i}+\frac{g_{i}G(N)}{N}\bigg)\\
				&-\frac{1}{2N}\ln x_{i}-\bigg(\frac{g_{i}G(N)}{N}-\frac{1}{2N}\bigg)\ln \frac{g_{i}G(N)}{N}+g_{i}\alpha\ln g_{i}\alpha+\frac{\omega(x_{i},N)}{N}\Bigg],\\
				\epsilon(N)=&
				\begin{cases}
					\frac{1}{N} \qquad &\text{if}\  \lim_{N\to\infty}\beta(N)N=0,\\
					|\beta(N)| \qquad &\text{otherwise}.
				\end{cases}
			\end{align*}
		\item Again, analogously to previous cases, by (\ref{entropy_properties_preposition_approximation_proof_approximation}) we have
			\begin{equation*}
				S(x,N)=G(N)f(x,N)+G(N)(\ln N+1)-\sum_{i=1}^{m}\bigg(g_{i}G(N)-\frac{1}{2}\bigg)\ln (g_{i}G(N))-\ln N,
			\end{equation*}
				where
			\begin{align*}
				f(x,N)=&\sum_{i=1}^{m}\Bigg[g_{i}\ln x_{i}+\bigg(\frac{x_{i}N}{G(N)}+g_{i}-\frac{1}{2G(N)}\bigg)\ln\bigg(1+\frac{g_{i}G(N)}{x_{i}N}\bigg)-g_{i}-\frac{\ln x_{i}}{G(N)}+\frac{\omega(x_{i},N)}{G(N)}\Bigg].
			\end{align*}
			Then
			\begin{equation*}
				f(x,N)=f(x)+\sigma(x,N)\epsilon(N),
			\end{equation*}
			where
			\begin{align*}
				f(x)=&\sum_{i=1}^{m}g_{i}\ln x_{i},\\
				\sigma(x,N)=&\frac{1}{\epsilon(N)}\sum_{i=1}^{m}\bigg[\bigg(\frac{x_{i}N}{G(N)}+g_{i}-\frac{1}{2G(N)}\bigg)\ln\bigg(1+\frac{g_{i}G(N)}{x_{i}N}\bigg)-g_{i}-\frac{\ln x_{i}}{G(N)}+\frac{\omega(x_{i},N)}{G(N)}\Bigg],\\
				\epsilon(N)=&
				\begin{cases}
					\frac{1}{G(N)} \qquad \text{if} \lim_{N\to\infty}\frac{G(N)^{2}}{N}=0,\\
					\frac{G(N)}{N} \qquad \text{otherwise}.
				\end{cases}
			\end{align*}
		\end{enumerate}
		By approximating the logarithm $\ln(1+x)=x+O(x^{2}), x\to 0$, one can verify that for all three cases, $\sigma(x,N)=O(1)$ as $N\to\infty$. Therefore the function $\sigma(x,N)$ is uniformly bounded on $\mathcal{A}^{(a)}_{E}$.\\
		\indent Moreover, for all three cases, it is clear that the functions $S(x,N)$, $f(x,N)$, $f(x)$ are at least three times differentiable.\\
		\indent Now, denote by $\Gamma(\lambda)$ Gamma function with $\lambda\in\mathbb{R}_{+}$. The logarithm of the Gamma function is a strictly convex function. This can be validated by checking the sign of the second derivative of $\ln\Gamma(\lambda)$. The derivative of Digamma function, $\psi(\lambda)=\Gamma'(\lambda)/\Gamma(\lambda)$, is equal to the second derivative of $\ln\Gamma$ and
		\begin{equation}
			\label{entropy_properties_preposition_approximation_proof_digamma}
			\frac{d^{2}\ln\Gamma(\lambda)}{d\lambda^{2}}=\frac{d\psi(\lambda)}{d\lambda}=\sum_{i=0}^{\infty}\frac{1}{(\lambda+i)^{2}}>0,
		\end{equation}
		which confirms strict convexity. For the series representation of the Digamma function see \cite{abramowitz_book} p. 259. \\
		\indent Functions $S(x,N)$, $f(x,N)$ and $f(x)$ are the sums of the $m$ components and each of these components depends only on single $x_{i}$'s, $i=1,\ldots,m$. Hence, cross derivatives are equal to zero and the Hessian matrices are diagonal matrices. We use (\ref{entropy_properties_preposition_approximation_proof_digamma}) to calculate the second derivatives of $S(x,N)$ given by (\ref{introduction_entropy}), that is
		\begin{align*}
			\frac{\partial^{2}S(x,N)}{\partial x_{i}^{2}}&=\frac{\partial^{2}}{\partial x_{i}^{2}}\ln\Gamma(x_{i}N+g_{i}G(N))-\frac{\partial^{2}}{\partial x_{i}^{2}}\ln\Gamma(x_{i}N-1)\\
			&=N\sum_{k=0}^{\infty}\bigg(\frac{1}{x_{i}N+g_{i}G(N)+k}-\frac{1}{x_{i}N-1+k}\bigg)\\
			&=-N\sum_{k=0}^{\infty}\frac{g_{i}G(N)+1}{(x_{i}N+g_{i}G(N)+k)(x_{i}N-1+k)},
		\end{align*}
		and for large enough $N$, these derivatives are negative. The function $f(x)$ has also negative second derivatives w.r.t $x_{i}$'s, that is, for each case of $G(N)$ and $x\in\mathcal{A}_{E}^{(a)}$ we have
		\begin{align*}
			1)\ &\frac{\partial f(x)}{\partial x^{2}}=-\frac{1}{x_{i}^{2}}<0,\\
			2)\ &\frac{\partial f(x)}{\partial x^{2}}=-\frac{g_{i}\alpha}{x_{i}(x_{i}+g_{i}\alpha)}<0,\\
			3)\ &\frac{\partial f(x)}{\partial x^{2}}=-\frac{g_{i}}{x_{i}^{2}}<0.\\
		\end{align*}
		Therefore, the Hessian matrices of $f(x)$ and, for large enough $N$, also of $S(x,N)$, are negative definite. This implies that those functions are strictly concave. Moreover, since function $f(x,N)$ is the only term in $S(x,N)$ that depends on $x$, hence for large enough $N$, $f(x,N)$ is also strictly concave.\\
		\indent Lastly, note that $\sigma(x,N)$ can be represented by $\sigma(x,N)=\sigma(x)+\omega'(N)O(1)$, where $\omega'(N)\to 0\ \text{as}\ N\to\infty$, and $\sigma(x)$ have the following form for each of case in (\ref{introduction_G_cases})
		\begin{align*}
			1)\ \sigma(x)=&
			\begin{cases}
				-m\ln\sqrt{2\pi}-\sum_{i=1}^{m}\frac{1}{2}\ln x_{i} \qquad &\text{if}\ \lim_{N\to\infty}\frac{N^{2}}{G(N)}=0,\\
				\sum_{i=1}^{m}\frac{2x_{i}^{2}}{g_{i}}\qquad &\text{otherwise},
			\end{cases}\\
			2)\ \sigma(x)=&
			\begin{cases}
				-m\ln\sqrt{2\pi}-\sum_{i=1}^{m}\frac{1}{2}\ln\frac{x_{i}+g_{i}\alpha}{x_{i}g_{i}\alpha} \qquad &\text{if}\ \lim_{N\to\infty}\beta(N)N=0,\\
				\sum_{i=1}^{m}g_{i}\big(2+\ln g_{i}(x_{i}+g_{i}\alpha)\big) \qquad &\text{otherwise},
			\end{cases}\\
			3)\ \sigma(x)=&
			\begin{cases}
				-m\ln\sqrt{2\pi}-\sum_{i=1}^{m}\ln x_{i} \qquad &\text{if}\ \lim_{N\to\infty}\frac{G(N)^{2}}{N}=0,\\
				\sum_{i=1}^{m}\frac{2x_{i}^{2}}{g_{i}}\qquad &\text{otherwise}.
			\end{cases}
		\end{align*}
		So, it is immediate that
		\begin{equation*}
			\sup_{x}\big|\sigma(x,N)-\sigma(x)\big|=\sup_{x}\big|O(1)\omega'(N)\big|\to 0,\ \text{as}\ N\to \infty,
		\end{equation*}
		hence $\sigma(x,N)$ is uniformly convergent.\\
		\indent Therefore, by theorem in \cite{lang_complex_analysis_book} p. 157, one has that also the derivatives of $\sigma(x,N)$ converges uniformly as $N\to\infty$, and so are uniformly bounded.
	\end{proof}

	\begin{proposition}
		\label{entropy_properties_preposition_optimization}
		The functions $f(x,N)$ and $f(x)$ from Proposition \ref{entropy_properties_preposition_approximation} have a unique maximum at $x^{*}(N)$ and $x^{*}$ respectively. Furthermore,
		\begin{enumerate}[1)]
		\item if $\varepsilon_{1}<E<\varepsilon^{T}g$, then $x^{*}(N),x^{*}$ are noncritical points on the boundary $\varepsilon^{T}x=E$ of the set $\mathcal{A}^{(a)}_{E}$. For each case of $G(N)$ given by (\ref{introduction_G_cases}), $x^{*}$ has the components
			\begin{align*}
				1)\ & x_{i}^{*}=\frac{g_{i}}{e^{\lambda\varepsilon_{i}+\nu}},\\
				2)\ & x_{i}^{*}=\frac{g_{i}\alpha}{e^{\lambda\varepsilon_{i}+\nu}-1},\\
				3)\ & x_{i}^{*}=\frac{g_{i}}{\lambda\varepsilon_{i}+\nu},
			\end{align*}
			for $i=1,\ldots,m$, where the parameters $\lambda>0$, $\nu$ are uniquely determined by the equations
			\begin{equation*}
				\sum_{i=1}^{m}x^{*}_{i}=1,\quad \varepsilon^{T}x^{*}=E,
			\end{equation*}
		\item if $E=\varepsilon^{T}g$ then the point $x^{*}$ is on the boundary. For the first case of $G(N)$, $x^{*}$ is a critical point of $f(x)$. For the second and third case $x^{*}$ is a noncritical point of $f(x)$,
		\item if $E>\varepsilon^{T}g$, then the points $x^{*}(N),x^{*}$ are in the interior of $\mathcal{A}^{(a)}_{E}$ and $x^{*}=g$.
		\end{enumerate}
	\end{proposition}

	\begin{proof}
		Let us consider optimization problem with strictly concave function $f(x)$ as an objective and the set $\mathbb{R}^{m}$ as the domain, that is
		\begin{align}
			\label{entropy_properties_preposition_optimization_proof_optimization_problem}
			\text{maximize}&\quad f(x),\notag\\
			\text{subject to}&\quad \mathbf{\varepsilon}^{T}x\leq E, \\
			&\quad \sum_{i=1}^{m}x_{i}=1, \notag\\
			&\quad x_{i}\geq 0,\ i=1,\ldots,m.\notag
		\end{align}
		By concavity, problem (\ref{entropy_properties_preposition_optimization_proof_optimization_problem}) can have at most one solution. For $E>\varepsilon_{1}$ there exists a strictly feasible point, hence the Slater's condition holds and the optimal vector exists. Therefore, the Karush-Kuhn-Tucker (KKT) conditions are necessary and sufficient for the existence and uniqueness of the optimal vector. For the details on the above Optimization Theory see \cite{convex_optimization}, Chapter 5.\\
		The KKT conditions for the problem (\ref{entropy_properties_preposition_optimization_proof_optimization_problem}) are the following
		\begin{align}
			\label{entropy_properties_preposition_optimization_proof_kkt}
			\varepsilon^{T}x^{*}&\leq E,\notag \\
			\sum_{i=1}^{m}x^{*}_{i}&=1,\notag\\
			\lambda&\geq 0,\\
			\lambda\Big( \varepsilon^{T}x^{*}-E\Big)&=0,\notag\\
			Df(x^{*})-\lambda\varepsilon-\nu&=0,\notag
		\end{align}
		where $D$ is a differential operator.\\
		The last condition, for each case of $G(N)$, yields
		\begin{align}
			\label{entropy_properties_optimization_last_condition_solution}
			1)\ & x_{i}^{*}=\frac{g_{i}}{e^{\lambda\varepsilon_{i}+\nu}},\notag\\
			2)\ & x_{i}^{*}=\frac{g_{i}c}{e^{\lambda\varepsilon_{i}+\nu}-1},\\
			3)\ & x_{i}^{*}=\frac{g_{i}}{\lambda\varepsilon_{i}+\nu}.\notag
		\end{align}
		Then, for the analysis of the conditions (\ref{entropy_properties_preposition_optimization_proof_kkt}) two cases can be distinguished, when $\lambda=0$ and when $\lambda>0$. In the first case we have
		\begin{align}	
			\label{entropy_properties_optimization_kkt_first_case_1}
			\varepsilon^{T}x^{*}&\leq E, \\
			\label{entropy_properties_optimization_kkt_first_case_2}
			\sum_{i=1}^{m}x^{*}_{i}&=1,\\
			\label{entropy_properties_optimization_kkt_first_case_3}
			\lambda&= 0.
		\end{align}
		From (\ref{entropy_properties_optimization_last_condition_solution}), (\ref{entropy_properties_optimization_kkt_first_case_2}) and (\ref{entropy_properties_optimization_kkt_first_case_3}) we get that the solution of the system is equal to  $x^{*}=g$. Moreover, by substituting that solution into (\ref{entropy_properties_optimization_kkt_first_case_1}) we obtain the condition for the parameter $E$, that is $E\geq\varepsilon^{T}g$. It can be easily inferred that when $E=\varepsilon^{T}g$, the point of maximum of $f$ is on the boundary of the set $\mathcal{A}^{(a)}_{E}$, which is the hyperplane $\varepsilon^{T}x=E$. When the parameter $E>\varepsilon^{T}g$, the point of maximum $x^{*}=g$ is in the interior of $\mathcal{A}^{(a)}_{E}$. Furthermore, for the first case of $G(N)$, since $\nu=0$, hence $Df(x^{*})=0$. Therefore the maximum is at a critical point of $f$. For the second and third case we have that $Df(x^{*})=\ln(1+c)$ and $Df(x^{*})=1$ respectively, hence the maximum of $f$ is at noncritical point.\\
		\indent When $\lambda>0$, we have the following system of equations
		\begin{align}
			\label{entropy_properties_optimization_kkt_second_case_1}
			\varepsilon^{T}x^{*}&=E,\\
			\label{entropy_properties_optimization_kkt_second_case_2}
			\sum_{i=1}^{m}x^{*}_{i}&=1,\\
			\label{entropy_properties_optimization_kkt_second_case_3}
			\lambda&> 0,
		\end{align}
		which has a unique solution for the remaining feasible values of $E$, that is, $\varepsilon_{1}<E<\varepsilon^{T}g$. Substitution of (\ref{entropy_properties_optimization_last_condition_solution}) into (\ref{entropy_properties_optimization_kkt_second_case_1}) and (\ref{entropy_properties_optimization_kkt_second_case_2}) yields system of equations from which we can obtain the parameters $\lambda$ and $\nu$. From (\ref{entropy_properties_optimization_kkt_second_case_1}) we infer that the maximum is on the boundary of $\mathcal{A}_{E}^{(a)}$. Furthermore, since $Df(x^{*})\neq 0$, the maximum is at noncritical point.\\
		\indent Finally, we consider the optimization problem (\ref{entropy_properties_preposition_optimization_proof_optimization_problem}) for the function $f(x,N)$ instead of $f(x)$. The solution of such problem is analogous. Although we are lacking of explicit expression for $x^{*}_{i}$'s, such as (\ref{entropy_properties_optimization_last_condition_solution}), the general structure of the solution is the same. Again, there are two cases, when $\lambda=0$ and when $\lambda>0$.\\
		\indent In the first case, i.e. $\lambda=0$, the maximum is attained at the point $x^{*}(N)=g(N)$ and is such that $g(N)\to g$ as $N\to\infty$. The maximum can be either on the boundary, where the parameter $E$ is equal to $E=\varepsilon^{T}g(N)$ or in the interior of $\mathcal{A}_{E}^{(a)}$ when $E>\varepsilon^{T}g(N)$.\\
		\indent In the second case of $\lambda$, the maximum $x^{*}(N)$ is attained for the remaining feasible values of the parameter $E$, that is $\varepsilon_{1}<E<\varepsilon^{T}g(N)$. For this situation, the point of maximum is on the boundary $E=\varepsilon^{T}x^{*}(N)$ and $x^{*}(N)\to x^{*}$ as $N\to\infty$.\\
	\end{proof}


\section{Limit theorems}
	\qquad Due to the constraint (\ref{particle_system_domain_constraint_1}), the last element of the vector $x$ is determined by the first $m-1$ elements, that is $x_{m}=\sum_{i=1}^{m-1}x_{i}$. So, within the set $\mathcal{A}_{E}$, the vector $x'=(x_{1},x_{2},\ldots,x_{m-1})$ uniquely determines the vector $x$. Therefore, we have that $S(x,N)=S(x',N)$, $f(x,N)=f(x',N)$ and consequently $f(x)=f(x')$ for any $x\in\mathcal{A}_{E}$. The situation is analogous for the random vector $X(N)$, and we conveniently denote $X'(N):=(X_{1}(N),X_{2}(N),\ldots,X_{m-1}(N))$. The pmf of $X'(N)$ is the following
	\begin{equation}
		\label{limit_theorems_pmf}
		P(X'(N)=x'):=\frac{e^{S(x',N)}}{\sum_{y'\in\mathcal{A}_{E}\cap L_{N}}e^{S(y',N)}}.
	\end{equation}
	Moreover, it can be easily verified that the entropy properties in Proposition \ref{entropy_properties_preposition_approximation} and \ref{entropy_properties_preposition_optimization} remain the same after the last element of $x$ is replaced by $x_{m}=1-\sum_{i=1}^{m-1}x_{i}$.\\
	\indent Now, let us introduce a preliminary result needed for the proofs of the limit theorems.
	\begin{proposition}
		\label{limit_theorems_preposition_pmf}
		For large enough $N$, the following approximation of the pmf (\ref{particle_system_pmf}) and (\ref{limit_theorems_pmf})  holds
		\begin{align*}
			&P(X(N)=x)=P(X'(N)=x')=\\
			&=
			\begin{cases}
				\frac{e^{h(N)f(x',N)}}{\sum_{\mathcal{A}^{(a)}_{E}\cap L_{N}}e^{h(N)f(y',N)}}\Big(1+O(1)e^{-h(N)\Delta}\Big),&\text{if}\ X'(N)\in\mathcal{A}^{(a)}_{E}\cap L_{N},\\
				O(1)e^{-h(N)\Delta},&\text{if}\ X'(N)\in\big(\mathcal{A}_{E}\backslash\mathcal{A}^{(a)}_{E}\big)\cap L_{N},
			\end{cases}
		\end{align*}
		where $\Delta>0$ is some constant, functions $h,f$ are given by Proposition \ref{entropy_properties_preposition_approximation}.
	\end{proposition}
	
	\begin{proof}
		We start by noticing that $P(X'(N)=x')=P(X(N)=x)$ because $x$ is uniquely determied by $x'$ for any $x\in\mathcal{A}_{E}$. Then, let us approximate the pmf (\ref{particle_system_pmf}) for $X(N)\in\mathcal{A}^{(a)}_{E}\cap L_{N}$.\\
		\indent First, we decompose the denominator
		\begin{equation}
			\label{limit_theorems_preposition_pmf_proof_decomposition}
			\sum_{\mathcal{A}_{E}\cap L_{N}}e^{S(x,N)}=\sum_{\mathcal{A}^{(a)}_{E}\cap L_{N}}e^{h(N)f(x,N)+C(N)}+\sum_{\big(\mathcal{A}_{E}\backslash\mathcal{A}^{(a)}_{E}\big)\cap L_{N}}e^{S(x,N)},
		\end{equation}
		where in the first sum we used the approximation of $S(x,N)$ from Proposition \ref{entropy_properties_preposition_approximation}.\\
		In Proposition \ref{entropy_properties_preposition_approximation} the function $f(x,N)$ is defined on $\mathcal{A}^{(a)}_{E}$. So, let us define $f(x,N)$ on $\mathcal{A}_{E}\backslash\mathcal{A}^{(a)}_{E}$ by
		\begin{equation}
			\label{limit_theorems_preposition_pmf_proof_extended_f}
			f(x,N):=\frac{1}{h(N)}(S(x,N)-C(N)),
		\end{equation}
		where $h(N)$ and $C(N)$ are some functions of $N$. Since $S(x,N)$ is strictly concave and differentiable w.r.t. $x$ on whole set $\mathcal{A}_{E}$, so $f(x,N)$ defined in (\ref{limit_theorems_preposition_pmf_proof_extended_f}) is also strictly concave and differentiable. Next, we apply Taylor's Theorem for $f(x,N)$ at $x\in\big(\mathcal{A}_{E}\backslash\mathcal{A}^{(a)}_{E}\big)\cap L_{N}$
		\begin{equation*}
			f(x,N)=f(x^{*},N)+Df(x_{\theta}(N),N)(x^{*}-x),
		\end{equation*}
		where $x_{\theta}(N)$ is some point between $x$ and $x^{*}$. Point $x^{*}$ is a point of maximum of $f$ defined in Proposition \ref{entropy_properties_preposition_optimization}. By the properties of the function $f(x,N)$ we have
		\begin{equation*}
			0<\Delta'<|Df(x_{\theta}(N),N)(x^{*}-x)|\leq\sup_{x\in\mathcal{A}_{E},N\geq N_{0}}\|Df(x,N)\|\sqrt{m}\leq\Delta,
		\end{equation*}
		since $|x^{*}-x|\leq\sqrt{m}$ as $x_{i}\in [0,1]$ for $i,\ldots,m$, and where $N_{0}\in\mathbb{Z}_{+}$ and $\Delta'$ such that $0<\Delta'<\Delta$ are some constants. Therefore
		\begin{equation*}
			\label{limit_theorems_preposition_pmf_proof_f_lower_bound}
			f(x,N)=f(x^{*},N)-|Df(x_{\theta}(N),N)(x^{*}-x)|\geq f(x^{*},N) - \Delta,
		\end{equation*}
		and
		\begin{equation}
			\label{limit_theorems_preposition_pmf_proof_f_upper_bound}
			f(x,N)\leq f(x^{*},N)-\Delta'.
		\end{equation}
		Then the second sum in (\ref{limit_theorems_preposition_pmf_proof_decomposition}) has a lower bound
		\begin{equation*}
			\sum_{\big(\mathcal{A}_{E}\backslash\mathcal{A}^{(a)}_{E}\big)\cap L_{N}}e^{S(x,N)}\geq e^{h(N)f(x^{*},N)+C(N)-h(N)\Delta}\sum_{\big(\mathcal{A}_{E}\backslash\mathcal{A}^{(a)}_{E}\big)\cap L_{N}}1.
		\end{equation*}
		Since the set $\mathcal{A}_{E}$ is bounded, the sum on the r.h.s. is also bounded, more precisely
		\begin{equation*}
			\sum_{\big(\mathcal{A}_{E}\backslash\mathcal{A}^{(a)}_{E}\big)\cap L_{N}}1\leq \sum_{\{x:0<x_{i}\leq 1; i=1,\ldots,m\}\cap L_{N}}=N^{m}.
		\end{equation*}
		Hence
		\begin{equation*}
			\sum_{\big(\mathcal{A}_{E}\backslash\mathcal{A}^{(a)}_{E}\big)\cap L_{N}}=O(1)N^{m}.
		\end{equation*}
		Combining above estimates yields
		\begin{equation}
			\label{limit_theorems_preposition_pmf_proof_denominator_estimate}
			\sum_{\mathcal{A}_{E}\cap L_{N}}e^{S(x,N)}\geq\sum_{\mathcal{A}^{(a)}_{E}\cap L_{N}}e^{h(N)f(x,N)+C(N)}+\omega(N)e^{C(N)},
		\end{equation}
		where 
		\begin{equation*}
			\omega(N):=O(1)N^{m}e^{h(N)f(x^{*},N)-h(N)\Delta}.
		\end{equation*}
		Hence, using (\ref{limit_theorems_preposition_pmf_proof_denominator_estimate}), the pmf for $X(N)\in\mathcal{A}^{(a)}_{E}\cap L_{N}$ can be estimated
		\begin{equation}
			\label{limit_theorems_preposition_pmf_proof_first_case_estimate}
			\frac{e^{S(x,N)}}{\sum_{\mathcal{A}_{E}\cap L_{N}}e^{S(y,N)}}\leq\frac{e^{h(N)f(x,N)}}{\sum_{\mathcal{A}^{(a)}_{E}\cap L_{N}}e^{h(N)f(y,N)}}\bigg(1-\frac{\omega(N)}{\sum_{\mathcal{A}^{(a)}_{E}\cap L_{N}}e^{h(N)f(y,N)}+\omega(N)}\bigg).
		\end{equation}
		Now, using Theorem \ref{appendix_theorem_sum_of_states_approximation_interior} and \ref{appendix_theorem_sum_of_states_approximation_boundary} from the Appendix we obtain
		\begin{equation}
			\label{limit_theorems_preposition_pmf_proof_sum_estimate}
			\sum_{\mathcal{A}^{(a)}_{E}\cap L_{N}}e^{h(N)f(x,N)}\geq O(1)N^{m}e^{h(N)f(x^{*},N)}h(N)^{-\frac{m+1}{2}},
		\end{equation}
		and substituting that into (\ref{limit_theorems_preposition_pmf_proof_first_case_estimate}) yields an estimate
		\begin{align*}
			&\Bigg|\frac{e^{S(x,N)}}{\sum_{\mathcal{A}_{E}\cap L_{N}}e^{S(y,N)}}-\frac{e^{h(N)f(x,N)}}{\sum_{\mathcal{A}^{(a)}_{E}\cap L_{N}}e^{h(N)f(y,N)}}\Bigg|\leq\frac{e^{h(N)f(x,N)}}{\sum_{\mathcal{A}^{(a)}_{E}\cap L_{N}}e^{h(N)f(y,N)}}\\
			&\times\bigg(\frac{\omega(N)}{\sum_{\mathcal{A}^{(a)}_{E}\cap L_{N}}e^{h(N)f(y,N)}+\omega(N)}\bigg)\leq\frac{\omega(N)h(N)^{\frac{m+1}{2}}}{\sum_{\mathcal{A}^{(a)}_{E}\cap L_{N}}e^{h(N)f(y,N)}}\\
			&\leq\frac{e^{h(N)f(x,N)}}{\sum_{\mathcal{A}^{(a)}_{E}\cap L_{N}}e^{h(N)f(y,N)}}O(1)e^{-h(N)\Delta''},
		\end{align*}
		where $\Delta''$ such that $0<\Delta''<\Delta'$, is some constant. Since $f(x,N)=f(x',N)$, hence we obtain the first part of the Proposition.\\
		\indent When $X(N)\in\big(\mathcal{A}_{E}\backslash\mathcal{A}^{(a)}_{E}\big)\cap L_{N}$ the pmf (\ref{particle_system_pmf}) can be approximated analogously. For the denominator we can use (\ref{limit_theorems_preposition_pmf_proof_denominator_estimate}) and (\ref{limit_theorems_preposition_pmf_proof_sum_estimate}). For the numerator, we use (\ref{limit_theorems_preposition_pmf_proof_f_upper_bound}) and obtain an estimate
		\begin{align*}
			\frac{e^{S(x,N)}}{\sum_{\mathcal{A}_{E}\cap L_{N}}e^{S(y,N)}}\leq&\frac{e^{h(N)f(x^{*},N)-h(N)\Delta'}}{\sum_{\mathcal{A}^{(a)}_{E}\cap L_{N}}e^{h(N)f(y,N)}+\omega (N)}\leq O(1)e^{-h(N)\Delta'}N^{-m}h(N)^{\frac{m+1}{2}}\\
			\leq& O(1)e^{-h(N)\Delta''}.
		\end{align*}
		Hence the Proposition is proved.
	\end{proof}

	Next, we provide the proofs of the limit theorems. The results also include estimates valid for sufficiently large $N$, where the parameter $\delta$ has values in the interval $\delta\in(0,\frac{1}{3(m+1)})$.

	\begin{theorem}[Weak law of large numbers]
		\label{limit_theorems_weak_law_of_large_numbers}
		As $N\to\infty$, the random vector $X(N)$ with the pmf (\ref{particle_system_pmf}) converges in distribution to $x^{*}$ and the following estimate of the mgf holds
		\begin{equation*}
			M_{X(N)}(\xi)=e^{\xi^{T}x^{*}}\bigg(1+\frac{O(1)}{h(N)^{1/2-3\delta}}+O(1)\frac{h(N)^{1/2+(m+1)\delta}}{N}+O(1)\epsilon(N)\bigg).
		\end{equation*}
	\end{theorem}

	\begin{remark}
		For this and the following limit theorems the convergence error term can be explicitly estimated using results in \cite{lapinski}.
	\end{remark}
	
	For the cases (\ref{introduction_G_cases}), let us assume $G(N)$ has the following properties
	\begin{align*}
		1)&\ \lim_{N\to\infty}\frac{N^{\frac{3}{2}}}{G(N)}=0,\\
		2)&\ \lim_{N\to\infty}\beta(N)\sqrt{N}=0,\\
		3)&\ \lim_{N\to\infty}\frac{G(N)^{\frac{3}{2}}}{N}=0.
	\end{align*}

	\begin{theorem}[Central limit theorem I]
		\label{limit_theorems_central_limit_theorem_1}
		Let $E>\varepsilon^{T}g$. For $X'(N)$ with distribution (\ref{limit_theorems_pmf}), the random vector $Z(N)=\sqrt{h(N)}(x'^{*}-X'(N))$ converges in distribution to $\mathcal{N}(0,D^{2}f(x'^{*})^{-1})$ and the following estimate of the mgf holds
		\begin{align*}
			M_{Z(N)}(\xi)=&\exp\bigg(\frac{1}{2}\xi^{T}D^{2}f(x'^{*})^{-1}\xi\bigg)\bigg(1+\frac{O(1)}{h(N)^{1/2-3\delta}}+O(1)\frac{h(N)^{1/2+(m+1)\delta}}{N}\\
			&+O(1)\epsilon(N)\sqrt{h(N)}\bigg), N\to\infty.
		\end{align*}
	\end{theorem}

	Here let us introduce notation $\xi_{y}=(\xi_{2},\ldots,\xi_{m-1})$, $Y=(X_{2}(N),\ldots,X_{m-1}(N))$ and $y^{*}=(x^{*}_{2},\ldots,x^{*}_{m-1})$.

	\begin{theorem}[Central limit theorem II]
		\label{limit_theorems_central_limit_theorem_2}
		Let $\varepsilon_{1}<E<\varepsilon^{T}g$ and consider the first or the second case of $G(N)$ in (\ref{introduction_G_cases}). For $X'(N)$ with distribution (\ref{limit_theorems_pmf}) there exists a subsequence of $N$ such that the random vector $$Z(N)=\big(N(x_{1}^{*}-X_{1}(N)),\sqrt{N}(y^{*}-Y(N))\big),$$ converges in distribution to a discrete distribution with the pmf
		\begin{equation*}
			P(Z_{1}(N)=i)=\exp\bigg(-\bigg|\frac{\partial f(x'^{*})}{\partial x_{1}}\bigg|i\bigg)\bigg(1-\exp\bigg(-\bigg|\frac{\partial f(x'^{*})}{\partial x_{1}}\bigg|\bigg)\bigg),
		\end{equation*}
		for $Z_{1}(N)$ and to $\mathcal{N}(0,D_{y}^{2}f(x'^{*})^{-1})$ for $\big(Z_{2}(N),\ldots,Z_{m-1}(N)\big)$.\\
		Furthermore, the following estimate of the mgf holds
		\begin{align*}
			M_{Z(N)}(\xi)&=\frac{1-\exp\big(-\big|\frac{\partial f(x'^{*})}{\partial x_{1}}\big|\big)}{1-\exp\big(-\big|\frac{\partial f(x'^{*})}{\partial x_{1}}\big|-\xi_{1}\big)}\exp\bigg(\frac{1}{2}\xi_{y}^{T}D_{y}^{2}f(x'^{*})^{-1}\xi_{y}\bigg)\times\\
			&\times\bigg(1+\frac{O(1)}{N^{1/2-3\delta}}+O(1)\frac{N^{1/2+(m+1)\delta}}{N}+\epsilon(N)\sqrt{N}O(1)\bigg), N\to\infty.
		\end{align*}
	\end{theorem}

	\begin{theorem}[Central limit theorem III]
		\label{limit_theorems_central_limit_theorem_3}
		Let $\varepsilon_{1}<E<\varepsilon^{T}g$ and consider the third case of $G(N)$ in (\ref{introduction_G_cases}). For $X'(N)$ with distribution (\ref{limit_theorems_pmf}) there exists a subsequence of $N$ such that the random vector $$Z(N)=\big(G(N)(x_{1}^{*}-X_{1}(N)),\sqrt{G(N)}(y^{*}-Y(N))\big),$$ converges in distribution to $Exp\big(\big|\frac{\partial f(x'^{*})}{\partial x_{1}}\big|\big)$ for $Z_{1}(N)$ and to $\mathcal{N}(0,D_{y}^{2}f(x'^{*})^{-1})$ for \\ $\big(Z_{2}(N),\ldots,Z_{m-1}(N)\big)$. Furthermore, the following estimate of the mgf holds
		\begin{align*}
			M_{Z(N)}(\xi)&=\frac{\big|\frac{\partial f(x'^{*})}{\partial x_{1}}\big|}{\big|-\big|\frac{\partial f(x'^{*})}{\partial x_{1}}\big|-\xi_{1}\big|}\exp\bigg(\frac{1}{2}\xi_{y}^{T}D_{y}^{2}f(x'^{*})^{-1}\xi_{y}\bigg)\\
			&\times\bigg(1+\frac{O(1)}{G(N)^{1/2-3\delta}}+O(1)\frac{G(N)^{1/2+(m+1)\delta}}{N}+O(1)\epsilon(N)\sqrt{G(N)}\bigg),N\to\infty.
		\end{align*}
	\end{theorem}
	
	\begin{remark}
		A particular subsequence of $N$, for which Theorem \ref{limit_theorems_central_limit_theorem_2} and \ref{limit_theorems_central_limit_theorem_3} are valid, can be found explicitly. For such subsequence the hyperplane of the boundary of $\mathcal{A}_{E}$ on which maximum of $f(x',N)$ is attained, coincide with the points on the boundary of the set $\mathcal{A}_{E}\cap L_{N}$. For the details how to find this subsequence see Remark \ref{appendix_remark_subsequence} in the Appendix.
	\end{remark}

	\begin{proof}[Proof of Theorem \ref{limit_theorems_weak_law_of_large_numbers}, \ref{limit_theorems_central_limit_theorem_1}, \ref{limit_theorems_central_limit_theorem_2}, \ref{limit_theorems_central_limit_theorem_3}]
		Let us consider the random vector $X'(N)$. First we approximate the pmf of $X'(N)$ given by (\ref{limit_theorems_pmf}) using Proposition \ref{limit_theorems_preposition_pmf}. Then, separately for each theorem, we use Proposition \ref{entropy_properties_preposition_approximation} and \ref{entropy_properties_preposition_optimization} to identify the properties of the functions $h(N)$ and $f(x',N)$ in the approximated pmf. \\
		\indent It is clear from Proposition \ref{limit_theorems_preposition_pmf} that the approximated pmf differs from the pmf (\ref{appendix_pmf_interior}) and (\ref{appendix_pmf_boundary}) only by the exponentially small term. Therefore, we can apply Theorems \ref{appendix_theorem_weak_law_of_large_numbers},  \ref{appendix_theorem_central_limit_theorem_1},  \ref{appendix_theorem_central_limit_theorem_2},  \ref{appendix_theorem_central_limit_theorem_3} from the Appendix for the random vector $X'(N)$ and to prove Theorems \ref{limit_theorems_weak_law_of_large_numbers}, \ref{limit_theorems_central_limit_theorem_1}, \ref{limit_theorems_central_limit_theorem_2}, \ref{limit_theorems_central_limit_theorem_3}, respectively. Theorems \ref{limit_theorems_central_limit_theorem_1}, \ref{limit_theorems_central_limit_theorem_2}, \ref{limit_theorems_central_limit_theorem_3} are proved immediately.\\
		In order to complete the proof of Theorem \ref{limit_theorems_weak_law_of_large_numbers} we need to transform the mgf of $X'(N)$ to the mgf of $X(N)$. We perform it in the following way
		\begin{align*}
			&e^{\xi_{m}}M_{X'(N)}(\hat{\xi}')=e^{\xi_{m}}M_{X'(N)}(\xi'-\xi_{m})=E\big[e^{(\xi'^{T}-\xi_{m})X'(N)+\xi_{m}}\big]=\\
			&=E\big[e^{\xi'^{T}X'(N)+\xi_{m}(1-X_{1}-\ldots-X_{m-1})}\big]=M_{X(N)}(\xi),
		\end{align*}
		where $\hat{\xi_{i}}=\xi_{i}-\xi_{m}$ and $\xi'=(\xi_{1},\xi_{2},\ldots,\xi_{m-1})$. We also need to transform the mgf of $x'^{*}$ to $x^{*}$
		\begin{equation*}
			e^{\xi_{m}}e^{\hat{\xi}'^{T}x'^{*}}=e^{\xi_{m}}e^{(\xi'-\xi_{m})^{T}x'^{*}}=e^{\xi'^{T}x'^{*}+\xi_{m}(1-x_{1}^{*}-\ldots-x_{m-1}^{*})}=e^{\xi'^{T}x'^{*}+\xi_{m}x_{m}^{*}}=e^{\xi^{T}x^{*}}.
		\end{equation*}
	Hence we proved the theorems.
	\end{proof}


\section{Overview of applications}
	\qquad The Zipf Law and it's extension, the Zipf-Mandelbrot law, are common laws in the Complexity Science. Behavior of a wide range of complex systems met in the real world can be modeled by these laws. For example, systems such as distribution of cities according to their population, word frequencies in the book, prices of car brands, and stock price changes on the market. For details on these examples see for example \cite{complexity_science_book, maslov_linguistics}. The power laws behavior of the above systems are obtained by fitting the system's data to the power law curves. The structure and properties of the systems has no direct input into modeling. Therefore, such an experimental approach is limited. It might work only for a particular, special case of the system and the data. For example, a small change of the system properties or limited data can invalidate model obtained in this way.\\
	\indent More robust and accurate modeling is possible if we understand the system, know it's structure and properties. Results in this paper might be helpful with such modeling. Particularly, Theorem \ref{limit_theorems_weak_law_of_large_numbers} provides a rigorous derivation of Zipf-Mandelbrot Law, based only on the underlying system properties. Therefore, such a result could be a tool to model the real system experiencing power law behavior and provide an explanation of such a behavior.\\
	\indent Aside of the derivation of the Zipf-Mandelbrot law, the main results of this paper shows that Maxwell-Boltzmann, Bose-Einstein and Zipf-Mandelbrot distributions can emerge from the same type of the system. Single system parameter determines which distribution occurs. Such unification should be a valuable result for the Theory of Complex Systems.\\
	\indent For example, let us consider the first two distributions obtained in Theorem \ref{limit_theorems_weak_law_of_large_numbers}, that is, Bose-Einsten and Maxwell-Boltzmann statistics. These are well known results of Classic and Quantum Statistical Mechanics. Most notably, these statistics are met in some gases under certain physical conditions, see for example \cite{ statistical_physics_pathria_book, statistical_physics_reif_book}. They also occur quite often for some complex systems. For example, Bose-Einstein statistics fits the distribution of the USA population income, see \cite{kusmartsev}. Another example is in \cite{wages_complexity_science}, there the study shows that UK population income very accurately fits Boltzmann distribution. Most importantly, the statistics of the system in \cite{kusmartsev} transforms from Bose-Einstain to Boltzman statistics when market condition changes. For the situation as in \cite{kusmartsev}, when two different statistics emerge from the same system with slightly different conditions or parameters, there might be a place to apply Theorem \ref{limit_theorems_weak_law_of_large_numbers}. The single system's parameter in the theorem, which determines whether it is Bose-Einstein or Boltzmann statistics, might be interpreted as market conditions of the system in \cite{kusmartsev}. Hypothetically, there might be another state of the system, that is, market conditions, for which the parameter is in different regime than for the first two cases. This another state could be the third distribution in Theorem \ref{limit_theorems_weak_law_of_large_numbers}, the Zipf-Mandelbrot law.\\

\begin{acknowledgments}
	Author dedicates special thanks to Professor Vassili Kolokoltsov from the University of Warwick for English language check, style suggestions and hint for the solution of the optimization problem.\\
	\indent Author would like to express gratitude to Professor Sergey Leble from Gda\'nsk University of Technology for the support in conducting this research.
\end{acknowledgments}
	

\appendix

\numberwithin{equation}{section}

\section{Appendix}
	\qquad Let us recall results from \cite{lapinski} with assumptions tailored for the results in this paper. We consider a bounded, open set $\mathcal{A}\subset\mathbb{R}^{m}$ and a lattice $L_{N}:=\{\frac{x}{N},x\in\mathbb{N}^{m}\}$ for $N\geq N_{0}$ with $N_{0}\in\mathbb{Z}_{+}$. We assume $0\in\mathbb{N}^{m}$. Then we introduce a function $f:\mathcal{A}\times\mathbb{Z}_{+} \to \mathbb{R}$ which have a unique maximum at $x^{*}(N)$ such that $\lim_{N\to\infty}x^{*}(N)=x^{*}$. Point $x^{*}$ is chosen to be the origin of our coordinate system. The derivatives of $f(\cdot,N)$ up to third order exists.  Moreover, we assume that 
	\begin{equation*}
		\Delta:=\inf_{N\geq N_{0}, x\in\mathcal{A}\backslash U}\{f(x^{*}(N),N)-f(x,N)\}>0,
	\end{equation*}
	where $U$ is some neighborhood of the origin.\\
	Additionally $f$ can be represented
	\begin{equation}
		\label{appendix_f_representation}
		f(x,N)=f(x)+\sigma(x,N)\epsilon(N),
	\end{equation}
	where $f(x), \sigma(x,N)$ are three times differentiable w.r.t. $x$, $\epsilon(N)>0$ for all $N\geq N_{0}$ and $\epsilon(N)\to 0$ as $N\to\infty$. We assume $\sigma(x,N)$ and its derivatives are uniformly bounded. Furthermore, we introduce a positive, increasing function $h:\mathbb{R}_{+}\to\mathbb{R}$ such that $\lim_{N\to\infty}\frac{h(N)}{N}=0$ or $h(N)=N$ and a differentiable function $g:\mathcal{A}\to\mathbb{R}$. Then
	\begin{enumerate}[(a)]
		\item for the sum
		\begin{equation}
			\label{appendix_sum_of_states_interior}
			\Sigma(N):=\sum_{\mathcal{A}\cap L_{N}}g(x)e^{h(N)f(x,N)},
		\end{equation}
		we assume $f(\cdot,N)$ and $f(\cdot)$ has a unique nondegenerate maximum in the interior of $\mathcal{A}$.
		\item for the sum
		\begin{equation}
			\label{appendix_sum_of_states_boundary}
			\Sigma(N):=\sum_{\mathcal{A}\cap L_{N}\cap\{x:x_{1}\geq 0\}}g(x)e^{h(N)f(x,N)},
		\end{equation}
		we assume $f(\cdot,N)$ and $f(\cdot)$ has a unique maximum on the boundary $\{x:x_{1}=0\}$. Additionally, $\frac{\partial f(x^{*}(N),N)}{\partial x_{1}}<0$, $\frac{\partial f(x^{*})}{\partial x_{1}}<0$. We also assume that on the boundary $\{ x:x_{1}=0\}$ and on every hyperplane parallel to that boundary, functions $f(\cdot,N)$ and $f(\cdot)$ have a unique nondegenerate maximum.
	\end{enumerate}
	
	\begin{remark}
		If the functions $f(x,N)$ and $f(x)$ are strictly concave then assumptions in the point (b) reduces to having maximum at a noncritical point in the interior of the boundary $\{x:x_{1}=0\}$.
	\end{remark}
	
	\begin{remark}
		\label{appendix_remark_subsequence}
		The situation when the boundary of the domain is $\{x:x_{1}=a\}$ with $a\in\mathbb{Q}_{+}$ can be reduced to the case of the boundary $\{x:x_{1}=0\}$, if we consider $N$ such that $Na\in\mathbb{Z}$. This is because for those values, the structure of the lattice $L_{N}$ is preserved after appropriate shift of the coordinate system.
	\end{remark}

	In the following theorems $\delta$ denote any number from the interval $(0,\frac{1}{3(m+1)})$
	\begin{theorem}
		\label{appendix_theorem_sum_of_states_approximation_interior}
		For the sum (\ref{appendix_sum_of_states_interior}), as $N\to\infty$, the following approximation holds
		\begin{align*}
			\sum_{\mathcal{A}\cap L_{N}}&g(x)e^{h(N)f(x,N)}=e^{h(N)f(x^{*}(N),N)}N^{m}\bigg(\frac{2\pi}{h(N)}\bigg)^{\frac{m}{2}}\\
			&\times\Bigg(\frac{g(x^{*}(N))}{\sqrt{|\det D^{2}f(x^{*}(N),N)|}}+\omega_{1}(N)\frac{1}{h(N)^{1/2-3\delta}}+\omega_{2}(N)G\frac{h(N)^{1/2+(m+1)\delta}}{N}\Bigg),
		\end{align*}
		where $\omega_{1}(N)=O(1)$ and $\omega_{2}(N)=O(1)$.
	\end{theorem}

	\begin{theorem}
		\label{appendix_theorem_sum_of_states_approximation_boundary}
		For the sum (\ref{appendix_sum_of_states_boundary}), as $N\to\infty$, the following approximation holds
		\begin{align*}
			&\sum_{\mathcal{A}\cap L_{N}\cap\{x_{1}:x_{1}\geq 0\}}g(x)e^{h(N)f(x,N)}=e^{N f(x^{*}(N),N)}N^{m-1}\bigg(\frac{2\pi}{h(N)}\bigg)^{\frac{m-1}{2}}\frac{1}{1-\exp\big(\frac{h(N)}{N}\frac{\partial f(x^{*}(N),N)}{\partial x_{1}}\big)}\\
			&\times\Bigg(\frac{g(x^{*}(N))}{\sqrt{\big|\det D_{y}^{2}f(x^{*}(N),N)\big|}}+\omega_{1}(N)\frac{1}{h(N)^{1/2-3\delta}}+\omega_{2}(N)\frac{h(N)^{1/2+(m+1)\delta}}{N}\Bigg),
		\end{align*}
		where $\omega_{1}(N)=O(1)$, $\omega_{2}(N)=O(1)$.
	\end{theorem}
	
	\begin{remark}
		For Theorem \ref{appendix_theorem_sum_of_states_approximation_boundary}, the situation when the boundary is an arbitrary hyperplane with rational coefficients can be reduced to the case with the boundary $\{ x:x_{1}= 0\}$. This is because after appropriate rotation of coordinate system, the structure of the lattice, which is essential for the application of the theorem is preserved. That is, all the points of the domain are on the equally spaced hyperplanes parallel to the boundary.
	\end{remark}

	Let $X(N)$ be a random vector with pmf defined using sums (\ref{appendix_sum_of_states_interior}) and (\ref{appendix_sum_of_states_boundary})
	\begin{align}
		\label{appendix_pmf_interior}
		(a)&\ P(X(N)=x):=\frac{e^{h(N)f(x,N)}}{\sum_{\mathcal{A}\cap L_{N}}e^{h(N)f(y,N)}},\\
		\label{appendix_pmf_boundary}
		(b)&\ P(X(N)=x):=\frac{e^{h(N)f(x,N)}}{\sum_{\mathcal{A}\cap L_{N}\cap\{x:x_{1}\geq 0\}}e^{h(N)f(y,N)}}.
	\end{align}

	\begin{theorem}[Weak law of large numbers]
		\label{appendix_theorem_weak_law_of_large_numbers}
		As $N\to \infty$, the random vector $X(N)$ converges in distribution to the constant $x^{*}$ and the following estimate of the mgf holds
		\begin{equation*}
			M_{X(N)}(\xi)=e^{\xi^{T}x^{*}}\bigg(1+\frac{O(1)}{h(N)^{1/2-3\delta}}+O(1)\frac{h(N)^{1/2+(m+1)\delta}}{N}+O(1)\epsilon(N)\bigg).
		\end{equation*}
	\end{theorem}

	For $\epsilon(N)=o\Big(\frac{1}{\sqrt{h(N)}}\Big)$, $N\to\infty$ we have the following results
	\begin{theorem}[Central limit theorem I]
		\label{appendix_theorem_central_limit_theorem_1}
		For $X(N)$ with distribution (\ref{appendix_pmf_interior}) the random vector $Z(N)=\sqrt{h(N)}(x^{*}-X(N))$ converges weakly to $\mathcal{N}(0,D^{2}f(x^{*})^{-1})$ and the following estimate of the mgf holds
		\begin{align*}
			M_{Z(N)}(\xi)=&\exp\bigg(\frac{1}{2}\xi^{T}D^{2}f(x^{*})^{-1}\xi\bigg)\bigg(1+\frac{O(1)}{h(N)^{1/2-3\delta}}+O(1)\frac{h(N)^{1/2+(m+1)\delta}}{N}\\
			&+O(1)\epsilon(N)\sqrt{h(N)}\bigg), N\to\infty.
		\end{align*}
	\end{theorem}

	Here let us introduce the notation $\xi_{y}=(\xi_{2},\ldots,\xi_{m})$, $Y=(X_{2}(N),\ldots,X_{m}(N))$ and $y^{*}=(x^{*}_{2},\ldots,x^{*}_{m})$.

	\begin{theorem}[Central limit theorem II]
		\label{appendix_theorem_central_limit_theorem_2}
		For $X(N)$ with distribution (\ref{appendix_pmf_boundary}) and assuming $h(N)=N$, the random vector $Z(N)=\big(N(x_{1}^{*}-X_{1}(N)),\sqrt{N}(y^{*}-Y(N))\big)$ converges weakly to a discrete distribution with pmf
		\begin{equation*}
			P(Z_{1}(N)=i)=\exp\bigg(-\frac{\partial f(x^{*})}{\partial x_{1}}i\bigg)\bigg(1-\exp\bigg(-\frac{\partial f(x^{*})}{\partial x_{1}}\bigg)\bigg),
		\end{equation*}
		for $Z_{1}(N)$ and to $\mathcal{N}(0,D_{y}^{2}f(x^{*})^{-1})$ for $\big(Z_{2}(N),\ldots,Z_{m}(N)\big)$. Furthermore, the following estimate of the mgf holds
		\begin{align*}
			M_{Z(N)}(\xi)&=\frac{1-\exp\big(-\frac{\partial f(x^{*})}{\partial x_{1}}\big)}{1-\exp\big(-\frac{\partial f(x^{*})}{\partial x_{1}}-\xi_{1}\big)}\exp\bigg(\frac{1}{2}\xi_{y}^{T}D_{y}^{2}f(x^{*})^{-1}\xi_{y}\bigg)\\
			&\times\bigg(1+\frac{O(1)}{N^{1/2-3\delta}}+O(1)\frac{N^{1/2+(m+1)\delta}}{N}+\epsilon(N)\sqrt{N}O(1)\bigg), N\to\infty.
		\end{align*}
	\end{theorem}

	\begin{theorem}[Central limit theorem III]
		\label{appendix_theorem_central_limit_theorem_3}
		For $X(N)$ with the distribution (\ref{appendix_pmf_boundary}) and assuming $\lim_{N\to\infty}\frac{h(N)}{N}=0$, the random vector $$Z(N)=\big(h(N)(x_{1}^{*}-X_{1}(N)),\sqrt{h(N)}(y^{*}-Y(N))\big),$$ converges weakly to $Exp\big|\frac{\partial f(x^{*})}{\partial x_{1}}\big|$ for $Z_{1}(N)$ and to $\mathcal{N}(0,D_{y}^{2}f(x^{*})^{-1})$ for \\
		$\big(Z_{2}(N),\ldots,Z_{m}(N)\big)$. Furthermore, the following estimate of the mgf holds
		\begin{align*}
			M_{Z(N)}(\xi)&=\frac{\big|\frac{\partial f(x^{*})}{\partial x_{1}}\big|}{\big|\frac{\partial f(x^{*})}{\partial x_{1}}-\xi_{1}\big|}\exp\bigg(\frac{1}{2}\xi_{y}^{T}D_{y}^{2}f(x^{*})^{-1}\xi_{y}\bigg)\\
			&\times\bigg(1+\frac{O(1)}{h(N)^{1/2-3\delta}}+O(1)\frac{h(N)^{1/2+(m+1)\delta}}{N}+O(1)\epsilon(N)\sqrt{h(N)}\bigg), N\to\infty.
		\end{align*}
	\end{theorem}
	

\bibliographystyle{plainnat}
\bibliography{biblio}

\end{document}